\documentclass[12pt]{amsart}
\usepackage{amssymb}
\usepackage{}
\usepackage{cases}

\usepackage{amsmath}
\usepackage{txfonts}
\theoremstyle{plain}
\newtheorem{Thm}{Theorem}

\newtheorem{Pro}[Thm]{Proposition}

\begin{document} 
\title[remarks on Bernstein theorems]
{remarks on volume growth and Bernstein theorems for translating solitons}

\author{Li Ma}

\address{Li Ma: Zhongyuan Institute of mathematics and Department of mathematics \\
Henan Normal university \\
Xinxiang, 453007 \\
China} \email{lma@tsinghua.edu.cn}

\dedicatory{}
\date{May 26th, 2015}

\begin{abstract}

In this paper, we study the properties of potential function of the translating soliton $M$ in $R^{n+1}$ and the volume growth of the intersection of Euclidean balls with $M$. We give a condition to obtain the Bernstein theorem for the translating solitons.
We also give an outline of a simple proof of the Bernstein theorem due to Bombieri-De Giorgi-Miranda.

{\textbf{Mathematics Subject Classification} (2000): 
53C44,35K58, 58J35}

{\textbf{Keywords}: Bernstein theorem, mean curvature flow, translating solitons,}
\end{abstract}

\thanks{The research is partially supported by the National Natural Science
Foundation of China No.11271111 and SRFDP 20090002110019.
}

\maketitle

\section{Introduction}\label{sect1}

In this note£¬we study the
Bernstein theorems for translating solitons. By definition, an n-dimensional hypersurface $M^n$ in the Euclidean space $R^{n+1}$ is called a translating soliton if it is a solution of mean curvature flow  $X_t=\bar{H}(X)=\Delta_MX$ (of hypersurfaces) in $R^{n+1}$ obtained by moving along the fixed direction $-a\in S^n\subset R^{n+1}$. Then the translating soliton equation under consideration is
$$ H=<\nu,a>$$
where $H$ is the mean curvature of $M$, $<.,.>$ is the inner product in $R^{n+1}$, and $\nu$ is the unit outward normal to $M$. Recall that the mean curvature vector field of $M$ in $R^{n+1}$ is $\bar{H}=-H\nu$ and $H=div(\nu)$ and $\Delta_M$ is the Laplacian operator on $M$ in the induced metric. We always assume that $M$ is properly immersed in $R^{n+1}$ and is complete with respect to the induced metric. It is well-known that the n-plane (with the vector $a\in R^{n+1}$ in it), the Grim Reaper $\Gamma\times R^{n-1}$ with $\Gamma=\{(x,y)\in R^2; y=-\log \cos x, x\in (-\pi/2,\pi/2)\}$, and the paraboloid soliton (also called "Bowl soliton") obtained by Altschuler-Wu \cite{AW} are translating solitons (see also \cite{Wa} and \cite{H} for uniqueness result).

We define the potential function $S(X)=<X,a>$ for the position vector $X \in M$ as in \cite{MV}. Below $\Delta_M$ is the Laplacian operator defined by induced Riemannian metric on $M$ and $\nabla S$ is the gradient of the function $S(X)$ on $M$. Then it is well-known that
$\bar{H}=\Delta_M X$.
As in \cite{MV} we have
$$
1=|a|^2=|a^T|^2+|a^\nu|^2=|\nabla S|^2+H^2
$$
and
$$
-\Delta_M S(X)=-<\Delta_M X,a>=H<\nu,a>=H^2, \ \ in \ \ M.
$$
 We denote by $B_R=B_R(0)=\{X\in R^{n+1}, |X|<R\}$ for $0\in M$. A difficult problem in the study of the translating soliton is to control the volume growth of $M\bigcap B_R$. Define
 $$\gamma=\overline{\lim}_{R\to \infty} \frac{1}{R^2}\log vol(M\bigcap B_R).$$
We remark that a related concept $\gamma_2$ has been introduced in L.Karp \cite{K}.

Then we have the following extension of a result in \cite{MV}.
\begin{Thm}\label{P1} Let $M$ be a translating soliton in $R^{n+1}$.

(1). There is an uniform positive constant $C>0$ such that
$$
\frac{1}{R}\int_{B_R(0)\bigcap M} |\nabla S|\geq C
$$
for all $R>0$.

(2). If $\gamma<\infty$, then either $\inf_MS=-\infty$ or $\sup_MS=\infty$.

(3). If $M$ is convex, then $\inf_M H=0$.

(4). If $M$ have uniformly bounded second fundamental form, then we have $\inf_MS=-\infty$.
\end{Thm}

We have the following consequence from Theorem \ref{P1} (2),i.e.,  if $vol(M\bigcap B_R)$ is of polynomial growth in $R>0$ and $\inf_MS>-\infty$ on the translating soliton $M$, we must have $\sup_MS=\infty$.
For a translating soliton defined by the graph of a smooth function on the plane we have a result below.

\begin{Pro}\label{P2} On a translating graphing soliton defined by the function $z=u(x)$ on some domain $\Omega\subset R^n$ with
$|Du|\leq C$ on $M$ for some uniform constant $C>0$, we have $\inf_M S=-\infty$. In other word, if the graph defined by the function $z=u(x)$ on $\Omega\subset R^n$ is a non-planar translating soliton with $\inf_M S>-\infty$, then we have $\sup_{M}|Du|= \infty$, where $Du=(u_x)$ and $|Du|=|u_x|$.
\end{Pro}

Recall that in the graph case, we have $\nu=(-Du,1)/\sqrt{1+|Du|^2}$ and
$$
H=-div(\frac{Du}{\sqrt{1+|Du|^2}}).
$$

By standard computation we know that the unit normal vector field $\nu$ on the hypersurface $M\subset R^{n+1}$ satisfies the equation
$$
\Delta_M \nu+|A|^2\nu =\nabla_MH.
$$
Sometimes we write by $\nabla=\nabla_M$. Since $H=<\nu,a>$ on the translating soliton $M$ and $\nabla_MH=\nabla_{e_j}H e_j$ for the local orthonormal moving frame $\{e_j\}$ on $M$, we have for any nontrivial constant vector $b\in R^{n+1}$,
$$
(\nabla_MH,b)=\nabla_{e_j}H (e_j,b)=(\nabla_{e_j}\nu,a)(e_j,b)=h_{ij}(e_i,a)(e_j,b)=II(a^T,b^T)
$$
where $II=(h_{ij}):=(<\nabla_i\nu,e_j>)$ is the second fundamental form of $M$ in the local frame $(e_j)$.
Using this and the log-trick (see \cite{CM}) we can get the following result for general solitons.

\begin{Pro}\label{P3} Assume that the two-dimensional soliton $M$ has the finite total curvature, i.e., $\int_M|A|^2<\infty$, where $A$ is the second fundamental form of $M$. Assume that there is a nontrivial constant vector $b\in R^{n+1}$ such that the function $f$ defined by $f(X)=(\nu,b)$ is non-negative on $M$, and $II(a^T,b^T)=h_{ij}(e_i,a)(e_j,b)\leq 0$ on $M$. Then $M$ is a hyperplane.
\end{Pro}

After presenting a proof of Proposition \ref{P3}, we outline a simpler proof of the famous Bernstein theorem due to Bombieri-De Giorgi-Miranda \cite{BDM}. Their results states as below.

\begin{Thm} If $M$ is a minimal surface defined by the entire function $u=u(x)$, where $x\in R^n$ with the assumption $\frac{1+|Du|^2}{(b', Du)^2+1}$ is uniformly bounded on $R^n$ for some vector $b'\in  R^n$, then $u$ is linear, that is to say, $M$ is plane.
\end{Thm}

Concerning with the area growth about the graphic translating solitons, we have the following general fact, which is more or less well-known to experts.

\begin{Thm}\label{P4} Let $M$ be an entire graphic surface defined by the function $u=u(x)$, $x\in R^n$. Assume that its mean curvature $H$ is non-positive. Then we have $vol(B_R\cap M)\leq CR^n$ for some uniform constant $C>0$. Furthermore, for such a graph being a translating soliton, we have either $\inf_MS=-\infty$ or $\sup_MS=\infty$.
\end{Thm}

It is quite possible to extend our result here to the expanding soliton defined by
$$
H=-<\nu,X>, \ \ on \ \ M.
$$
In fact, we can follow the argument of Theorem 2.1 and Cor. 2.3 in M.Anderson \cite{A} to prove the following
result.
\begin{Pro}\label{P5} Let $M^n\to R^N$
be a complete expanding soliton immersion. There exists a constant $\epsilon_0=\epsilon_0(n,N)>0$ such that if
$$
\int_M |A|^ndv_g\leq \epsilon_0,
$$
then $M$ is an affine plane in $R^N$.
\end{Pro}
Since the argument of Proposition \ref{P5} just goes as in the same way as in Cor. 2.3 \cite{A} via a use of well-known compactness result and an epsilon regularity for expanding solitons, we omit the full detail. This fact was observed by us with Dr.Anqiang Zhu. We expect that with the assumption of finite total curvature, i.e., $\int_M |A|^ndv_g<\infty$, the structure
of the expanding soliton $(M,g)$ should be as good as minimal immersions.

\section{arguments of Propositions and Theorems \ref{P1} and  \ref{P4}}\label{sect2}
We may assume $0\in M$.

\emph{Proof of Theorem \ref{P1}}:

(1). We argue by contradiction. That is,
 there exists a positive number sequence $\{R_j\}$ with $R_j\to \infty$ such that
$\frac{1}{R}\int_{B_R(0)\bigcap M} |\nabla S|\to 0$ for $R=R_j\to\infty$.

Let $\phi_R(X)$ be a cut-off function defined on $B_R(0)$ and $\phi_R(X)=1$ on $B_{R/2}(0)$. Then we have
$$
\int_{B_{R/2}\bigcap M}H^2\leq \int H^2\phi_R^2=\int \phi_R<\nabla S(X),\nabla \phi_R(X)>
$$
which is bounded by
$$
\frac{1}{R}\int_{B_R(0)\bigcap M} |\nabla S|\to 0
$$
as $R\to \infty$. Then $H=0$ and $<\nu,a>=0$ on $M$. This implies that $|\nabla S(X)|=1$ on $M$, which is impossible by the result in \cite{MV}.

(2).  We may use Theorem 2.3 in \cite{K}. Note that the geodesic ball of radius $R$ at center $0$ is always contained in $B_R(0)\bigcap M$, our condition implies the condition $\gamma_2\leq \gamma<\infty$ in there. Here $\gamma_2$ is defined in p.450 in \cite{K}. Note that
$$
-\Delta_MS=H^2=1-|\nabla S|^2.
$$
We let $w=e^{-S}$. Then we have $\Delta_M w=w$. Assume that $\inf_M S>-\infty$. Then $\sup_M w<\infty$. Then by Theorem 2.3 in \cite{K}, we have
$\inf_M w=\inf_M\Delta_Mw\leq 0$. Remember that $\inf_M w\geq 0$. Then $\inf_Mw=0$, which is equivalent to $\sup_M S=\infty$.

(3). Let $u(x)=-S(x)$ on $M$. Then $\Delta_Mu=H^2$ and $|\nabla u|\leq 1$ in $M$. By the convexity of the hypersurface $M$ in $R^{n+1}$,
we know that $RicM\geq 0$ on $M$. Hence we can apply Cor. 2.2.2 in \cite{K} to conclude that
$\inf_M H^2=\inf_M(\Delta u)\leq 0$.

(4).  Assume now that $d:=\inf_M S>-\infty$. Then $\sup_M w=exp(-d)$. By our assumption, we know that the Ricci curvature of $M$ is uniformly bounded below so that the Omori-Yau maximum principle holds true on $M$. That is there is a point sequence $(x_j)$ in $M$ such that
$$
w(x_j)\to exp(-d), \ \  |\nabla w(x_j)|\leq \frac{1}{j}, \ \ \text{and} \ \  \Delta_Mw(x_j)\leq \frac{1}{j},
$$
as $j\to\infty$.
 By the equation $\Delta_M w=w$ we know that $\frac{1}{j}\geq \Delta_Mw(x_j)=w(x_j)\to exp(-d)>0$, which is impossible.
\qed

We remark that the proof of Proposition 1(4) can also be used to prove the conclusion \cite{MV} that\emph{ if the two dimensional translating soliton $M^2$ is conformal to the plane, then $\inf_M S=-\infty$}.
In fact, from $-\Delta_M (exp(-d)-w)=w>$ on $M$, we know that $exp(-d)-w>0$ is a positive superharmonic function on $R^2$. By the Liouville theorem we know that $w$ is a constant function so that $S$ is a constant function on $M$, which is impossible.

\emph{Proof of Proposition \ref{P2}}:  We argue by contradiction, i.e., $\inf_MS>-\infty$. We may assume that $H\not=0$ somewhere in $M$. For otherwise, it follows from the well-known Bernstein theorem \cite{CM} that $M$ is a hyperlane and then $\inf_M S=-\infty$. As we have assumed, $M$ is a graph defined the function $z=u(x)$ on the domain $\Omega\subset R^n$. Recall that
$$
\Delta_Mw=w, \ \ in \ \ M,
$$
where $w=exp(-S)$.
Since $Vol(B_R\cap M)=\int_{B_R\cap M} \sqrt{1+|Du|^2}\leq CR^n$. We remark that by Proposition \ref{P1} (2), we know
$\sup_MS=\infty$.  Let $c=\inf_M S$ and $B=e^{-c}=\sup_M w$.
Let $B(R)$ be the geodesic ball of radius $R$ with center $0$.
Let $f(R)=\int_{B(R)} w$. Then $f'(R)=\int_{\partial B(R)} w $. Note that $\nabla w=-w\nabla S$.
Integrating the equation $\Delta_Mw =w$ over $B(R)$ we obtain
$$
f(R)=\int_{B(R)}\Delta w=-\int_{\partial B(R)} <\nabla w, \partial/\partial_R>\leq \int_{\partial B(R)}w=f'(R),
$$
which implies that for any $R>R_0>0$,
\begin{equation}\label{vol-1}
f(R)\geq f(R_0)exp(R-R_0).
\end{equation}
However, by our assumption,
$$f(R)\leq B vol(B(R))\leq B vol(M\bigcap B_R)\leq CBR^n.$$ This gives a contradiction to (\ref{vol-1}) when $R$ large.

 \qed

\emph{Proof of Proposition \ref{P3}}: By the strong maximum principle we may assume $(\nu,b)>0$ on $M$. Since $M$ has the finite total curvature, we have the quadratic area growth, which is well-known to experts and may be due to many people, see the works of Cohn-Vossen, M.Anderson, and B.White \cite{W} \cite{CM}. 
Recall that
$$
\Delta_M \nu+|A|^2\nu=\nabla H, \ \ on \ \ M.
$$
Note that
$$
<\nabla H,b>=II<a^T,b^T>.
$$
By $II(a^T,b^T)\leq 0$, we get the differential inequality
$$
\Delta_M (\nu,b)+|A|^2(\nu,b)\leq 0, \ \ on \ \ M.
$$
By this (see \cite{CM} for detail) we can get the stability condition
$$
\int_M |A|^2 \eta^2\leq \int_M |\nabla_M\eta|^2
$$
for any $\eta\in C_0^1(M)$. Then using the log-trick \cite{CM}  we can get that $|A|=0$, which implies that $M$ is a plane. \qed

We now point out that we can use the idea above to give a simple proof of the Bernstein theorem due to Bombieri-De Giorgi-Miranda \cite{BDM}, which says that
if $M$ is a minimal surface defined by the entire function $u=u(x)$ $,x\in R^n$ with $\frac{1+|Du|^2}{(b', Du)^2+1}$ is uniformly bounded on $M$ for some vector $b'\in  R^n$, then $u$ is linear, that is to say, $M$ is plane.

\textbf{Here is an outline of the simpler proof}. Recall the formula for the smooth functions $f$, $h>0$ on $M$, and for $\Delta=\Delta_M$,
$$
\Delta{\frac f h}=\frac{\Delta f}{h}-\frac{f\Delta h}{h^2}-\frac{2}{h}\nabla_{\nabla h}\frac {f}{h}.
$$
Write by $\nu=(\nu_1,...,\nu_{n+1})$, where $\nu_j=-\frac{D_ju}{\sqrt{1+|Du|^2}}$ for $j=1,...,n$ and $\nu_{n+1}=\frac{1}{\sqrt{1+|Du|^2}}$.
Let $f(X)=(b,\nu)$ where $b=(b', b_{n+1})\in R^{n+1}$ is a fixed vector and let $h=\nu_{n+1}=\frac{1}{v}$ where $v=\sqrt{1+|Du|^2}$. Then we have
$$
\Delta(f/h)=(\nabla (f/h), \nabla \log v^2).
$$
Then we have
$$
div_M (v^{-2}\nabla (b',Du))=0,
$$
which can be written as
$$
div_M [\frac{(b', Du)^2+1}{1+|Du|^2}\nabla g]=0
$$
for $g=arctg (b',Du)$. Assume that $\frac{1+|Du|^2}{(b', Du)^2+1}$ is bounded by some constant on $M$, we then can do the Moser iteration \cite{BG} to conclude the Harnack inequality for $g$ and hence $g$ is constant, which in turn implies that $|Du|$ is bounded. We then conclude by Moser theorem that $u$ is linear.

\emph{Theorem \ref{P4}} can be proved by using the formula
$$
H=-div(\frac{Du}{\sqrt{1+|Du|^2}})\leq 0.
$$
Note that for any $c\in R$, $u+c$ has same mean curvature as $u$.
With loss of generality we may assume that $u(0)=0$ and we can using the truncation function
$u_R+R\geq 0$ by the truncation process that $u$ does not when $|u|\leq R$, $u_R=-R $ for $u\leq -R$, and $u=R$ for $\geq R$.  Multiplying both sides of the mean curvature formula
\begin{equation}\label{MS}
H=-div(\frac{Du}{\sqrt{1+|Du|^2}})\leq 0
\end{equation}
by the non-negative function $u_R+R$ and integrating part over $D:=D_R=\{x\in R^n; |x|<R\}$, we have
$$
\int_{D}\frac{|Du|^2}{\sqrt{1+|Du|^2}}\leq \int_{\partial D}\frac{|u_R+R||Du|}{\sqrt{1+|Du|^2}},
$$
which implies that
$$
vol(B_R\cap M)\leq\int_D\sqrt{1+|Du|^2}\leq CR^n.
$$

Applying Proposition \ref{P1} (2) we conclude that either $\inf_MS=-\infty$ or $\sup_MS=\infty$.
 This completes the proof of Theorem \ref{P4}.

\end{document}